\documentclass[12pt]{amsart}
\usepackage[mathscr]{eucal}
\usepackage{hyperref}

\newtheorem{theorem}{Theorem}[section]
\newtheorem{lemma}[theorem]{Lemma}
\newtheorem{proposition}[theorem]{Proposition}
\newtheorem{corollary}[theorem]{Corollary}
\newtheorem*{proposition 3.4}{Proposition \ref{Kroneckeropen}}

\theoremstyle{definition}
\newtheorem{definition}[theorem]{Definition}
\newtheorem*{notation}{Notation}
\newtheorem*{claim}{Claim}

\newtheorem{question}[theorem]{Question}

\theoremstyle{remark}
\newtheorem*{remark}{Remark}

\begin{document}
\title[Sumsets of dense sets and sparse sets]{Sumsets of dense sets and sparse sets}
\author[John T. Griesmer]{John T. Griesmer}
\address{Department of Mathematics\\
University of British Columbia\\
1984 Mathematics Road\\
Vancouver, BC V6T 1Z2\\
Canada}
\email{griesmer@math.ubc.ca}
\date{\today}

\begin{abstract}
  R.~Jin showed that whenever $A$ and $B$ are sets of integers having positive upper Banach density, the sumset $A+B:=\{a+b:a\in A, b\in B\}$ is piecewise syndetic.  This result was strengthened by Bergelson, Furstenberg, and Weiss to conclude that $A+B$ must be piecewise Bohr.  We generalize the latter result to cases where $A$ has Banach density $0,$ giving a new proof of the previous results in the process.
\end{abstract}

\maketitle

\setcounter{tocdepth}{1}%%Only sections appear in Contents
\tableofcontents

\section{Introduction}\label{Introduction}

\subsection{Large sets of integers and sumsets}

The well-known Steinhaus lemma says that whenever $A$ and $B$ are subsets of a locally compact group $G$ both having positive Haar measure, the product set $\{ab:a\in A, b\in B\}$ contains a nonempty open subset of $G.$  R.~Jin found an interesting analogue of the Steinhaus lemma for subsets of $\mathbb Z.$  To state his result, we need the notions of upper Banach density and piecewise syndeticity.

\begin{notation}  If $A$ and $B$ are subsets of an abelian group $G,$ and $c\in G,$ we write $A+c$ for $\{a+c:a\in A\},$ and $A+B$ for $\{a+b:a\in A,b\in B\}.$
\end{notation}

\begin{definition}  Let $A\subset \mathbb Z.$  The \textit{upper Banach density} of $A$ is the number
$$
d^*(A):=\lim_{M\to\infty}\sup_{N\in \mathbb Z} \frac{|A\cap [N,M]|}{N-M+1},
$$
and the \textit{upper density} of $A$ is the number
$$
\bar d(A):=\limsup_{N\to \infty} \frac{|A\cap [1,N]|}{N},
$$
where $|S|$ denotes the cardinality of the set $S.$

We also use the following standard terminology;  see \cite{BFW} for elaboration.

\begin{enumerate}

\item[$\bullet$] A set $S\subset \mathbb Z$ is \textit{syndetic} if there is a finite set $F\subset \mathbb Z$ such that $\mathbb Z=S+F.$  Equivalently, $S$ is syndetic if there exists $L$ such that $S\cap [M,N]$ is nonempty whenever $N-M>L.$

\item[$\bullet$]  A set $S\subset \mathbb Z$ is \textit{thick} if for every $L,$ $S$ contains an interval of length $L.$

\item[$\bullet$] A set $S\subset \mathbb Z$ is \textit{piecewise syndetic} if it is the intersection of a syndetic set with a thick set.
\end{enumerate}
Note that piecewise syndetic sets are always nonempty.
\end{definition}

\begin{theorem}[\cite{Jin01}, Corollary 3]\label{Jin1}

If $A,B\subset \mathbb Z$ with $d^*(A)>0, d^*(B)>0,$ then $A+B$ is piecewise syndetic.
\end{theorem}
Jin deduced this result, and the Steinhaus lemma (for $\mathbb R$), from a theorem in nonstandard analysis.

\medskip

Theorem \ref{Jin1} was strengthened in \cite{BFW}, where the conclusion that $A+B$ is piecewise syndetic was replaced by the conclusion that $A+B$ is piecewise Bohr.  We will discuss Bohr and piecewise Bohr sets in Section \ref{bs}. Briefly, a set $S\subset\mathbb Z$ is \textit{Bohr} if there is a trigonometric polynomial $P(n)=\sum_{j=1}^k c_je^{i\lambda_j n},$ $\lambda_j\in \mathbb R,$ such that $S\supset \{n: \operatorname{Re} P(n)>0\}\neq \emptyset,$ while a set is \textit{piecewise Bohr} if it is the intersection of a thick set and a Bohr set.

\begin{theorem}[\cite{BFW}, Theorem I]\label{BFWthm}  If $A,B\subset \mathbb Z$ with $d^*(A)>0, d^*(B)>0,$ then $A+B$ is piecewise Bohr.
\end{theorem}

Our main theorem generalizes Theorem \ref{BFWthm} by weakening the hypothesis that $d^*(A)>0,$ replacing upper Banach density by a more general notion of density satisfying an equidistribution condition.  Given a sequence $(\nu_j)_{j\in \mathbb N}$ of probability measures  on $\mathbb Z,$ one can consider the upper density $d_\nu$ with respect to $(\nu_j)_{j\in \mathbb N}$:  for $A\subset \mathbb Z$ define
$$
d_\nu(A):=\limsup_{j\to \infty} \nu_j(A).
$$
We say that a sequence  $(\nu_j)_{j\in \mathbb N}$ of probability measures on $\mathbb Z$ is an \textit{equidistributed averaging sequence} if for all $\theta\in (0,2\pi),$ we have
$$
\lim_{j\to \infty} \int e^{in\theta}\, d\nu_j(n)=0.
$$

With this definition we can state our main theorem.

\begin{theorem}\label{eastheorem} Let $(\nu_j)_{j\in \mathbb N}$ be an equidistributed averaging sequence, and let $A, B\subset \mathbb Z$ with $d^*(B)>0.$  Then the following implications hold.
\begin{enumerate} \item[1.]  If $d_\nu(A)>0,$ then $A+B$ is piecewise Bohr.
\item[2.]  If $d_\nu(A)=1,$ then $A+B$ is thick.
\end{enumerate}
\end{theorem}

\begin{remark} The hypothesis that $(\nu_j)_{j\in \mathbb N}$ is equidistributed is analogous to the hypothesis of $\alpha$-uniformity, or pseudorandomness, exploited in \cite{Gre2}, \cite{GT}, \cite{HL}, \cite{TV}, and other work on additive combinatorics.  Lemma 4.13 of \cite{TV} relates pseudorandomness to sumsets $A_1+A_2+\cdots +A_k$ where $k\geq 3,$ and is similar in spirit to Theorem \ref{eastheorem}.

Equidistributed averaging sequences are heavily studied in ergodic theory; see Section \ref{ergodicsection}.

The measures $\nu_j$ we have in mind take the form $\nu_j=$ normalized counting measure on $A_j,$ where $A_j\subset \mathbb Z$ is a finite set, so the averages $\int e^{in\theta}\, d\nu_j(n)$ can be written as $\frac{1}{|A_j|}\sum_{n\in A_j} e^{in\theta}.$  We use the notation of measures and integrals, rather than sets and averages, because our proofs seem more motivated with this notation.
\end{remark}

For certain choices of $\nu_j,$ Theorem \ref{eastheorem} produces interesting corollaries.   For instance,  $\nu_j=$ normalized counting measure on an interval $I_j,$ with $|I_j|\to \infty$, defines an equidistributed averaging sequence, and so we recover Theorems \ref{Jin1} and \ref{BFWthm}.  As we shall see in Section \ref{preliminaries}, a family of equidistributed averaging sequences is given by $\nu_j=$ normalized counting measure on $\{\lfloor 1^\alpha \rfloor , \lfloor 2^{\alpha} \rfloor, \dots, \lfloor j^\alpha \rfloor\},$ where $\alpha>0$ is not an integer, so we obtain the following corollary.
\begin{corollary}\label{powers}
If $\alpha>0$ is not an integer and $d^*(B)>0,$ then
$$
\{\lfloor n^{\alpha} \rfloor+b: n\in \mathbb N, b\in B\}
$$
is thick.  Furthermore, if $C\subset \mathbb N$ with  $\bar d(C)>0,$ and $d^*(B)>0,$ then
$$\{\lfloor c^{\alpha} \rfloor+b : c\in C, b\in B\}$$
is piecewise Bohr.
\end{corollary}
For other examples of equidistributed averaging sequences, and hence more results like the above corollary, one may consult \cite{BKQW} and its bibliography.

\subsection{Examples}\label{introexamples}  In contrast with Corollary \ref{powers}, we will derive the following examples from Proposition \ref{examples}.

\medskip

\begin{enumerate}
\item[$\bullet$] If $k\in \mathbb N$ and $k\geq 2,$ then for all $\varepsilon>0,$ there exists $B\subset \mathbb Z$ with $d^*(B)>1-\varepsilon$ such that $\{n^k+b:n\in \mathbb N, b\in B\}$ is not piecewise syndetic.
\medskip
\item[$\bullet$]  For all $\varepsilon>0,$ there exists $B\subset \mathbb Z$ with $d^*(B)>1-\varepsilon$ such that \\ $\{p+b: p \text{ is prime}, b\in B\}$ is not piecewise syndetic.
\end{enumerate}

\subsection{Arithmetic progressions.}  With little extra effort, the proof of Theorem \ref{eastheorem} establishes the following fact about the density of arithmetic progressions in sumsets.

\begin{theorem}\label{APs}  Let $(\nu_j)_{j\in \mathbb N}$ be an equidistributed averaging sequence, let $A,B\subset \mathbb Z$ with $d_\nu(A)>0, d^*(B)>0.$  Then for all $k\in \mathbb N$ and all $\varepsilon>0,$ there exists $d\in \mathbb N$ such that
\begin{align}\label{AP<}
d^*\Bigl(\bigcap_{l=0}^k \bigl((A+B)-ld\bigr)\Bigr)>\max\{d_\nu(A),d^*(B)\}-\varepsilon.
\end{align}
In fact, there is a Bohr set of $d$ satisfying (\ref{AP<}).
\end{theorem}
The conclusion says that there are many $d$ such that the set of $c$ with $\{c,c+d,c+2d,\dots, c+kd\}\subset A+B$ has density at least $\max\{d_\nu(A),d^*(B)\}-\varepsilon.$  This exhibits a sharp distinction between sumsets of the kind we consider and arbitrary sets of positive upper Banach density: \cite{BHK} features a construction of sets $B$  with $d^*(B)>0$ having much less than the expected density ($d^*(B)^5$) of $5$-term arithmetic progressions $\{c,c+d,c+2d,c+3d,c+4d\}$ for every nonzero $d.$    We will elaborate on this topic in Section \ref{bs}; for now we remark that in the case $d^*(A), d^*(B)>0,$ Theorem \ref{APs} is implicit in the proof of \cite{BFW}, Theorem I.

\begin{remark}  In \cite{BBF}, Theorems \ref{Jin1} and \ref{BFWthm} are generalized to the setting where $\mathbb Z$ is replaced by an arbitrary countable amenable group, answering a question in \cite{JK}.  Also, \cite{BBF} shows that the conclusion of Theorem \ref{BFWthm} cannot be qualitatively strengthened, as every piecewise Bohr set contains a set of the form $A+B$, where $d^*(A),d^*(B)>0.$

In his Ph.D. thesis (\cite{Gri}), the author proved a version of Theorem \ref{eastheorem} in the setting where $\mathbb Z$ is replaced by a countable amenable group.  This article confines itself to $\mathbb Z,$ because the methods of \cite{Gri} are very similar to those here, while the notation there is more cumbersome.  Furthermore, the examples and questions we formulate are most easily understood in the integer setting.
\end{remark}

\begin{remark} All of the aforementioned results are instances of the general phenomenon that whevever $A$ and $B$ are subsets of some (perhaps nonabelian) group, the sumset $A+B$ (product set $A\cdot B$) tends to be more structured than $A$ and $B$.  Quantitative examples of this phenomenon are discussed in \cite{B}, \cite{CS}, \cite{Gre1}, \cite{HL}, and \cite{San}, which deal with finding long arithmetic progressions in sumsets, while \cite{CH} considers sumsets $A+A$ where $A$ is a set of primes.  The book \cite{TV} considers this and many related issues, and has an extensive bibliography.

\end{remark}

\subsection{Outline of the article}

Theorems \ref{eastheorem} and \ref{APs} will be deduced from Theorem \ref{eastheorem'}, which is a slight refinement of Theorem \ref{eastheorem}.  Theorem \ref{eastheorem'} will be deduced from Proposition \ref{Kroneckeropen}, the latter being an ergodic theoretic analogue of Theorem \ref{eastheorem'}.  We invoke the Steinhaus lemma at a crucial step in the proof of Proposition \ref{Kroneckeropen};  this seems to strengthen the analogy between sumsets in $\mathbb Z$ and sumsets in general locally compact groups.

We will see that the problem of describing sumsets $A+B,$ where $d^*(B)>0,$ corresponds loosely to the problem of identifying limits of averages of the form $\frac{1}{|A_j|}\sum_{a\in A_j} f\circ T^{-a},$ where $(X,\mathscr X,\mu,T)$ is a measure preserving system, $f:X\to [0,1],$ and the $A_j$ are finite subsets of $A.$  Under the hypothesis that $(\nu_j)_{j\in \mathbb N}$ is an equidistributed averaging sequence and $d_\nu(A)>0,$ we obtain a fairly precise description of such limits in the proof of Proposition \ref{Kroneckeropen}.

%remark about embeddings?

In the next section we summarize the definitions and tools we need from ergodic theory, pausing in Subsection \ref{bs} to delineate some differences between piecewise Bohr sets and arbitrary sets of positive density. In Section \ref{ergodicschemes} we prove Theorem \ref{eastheorem}, and in Section \ref{examplessection} we provide the examples promised in Section \ref{introexamples}.  In Section \ref{questions} we ask some natural questions raised by Theorem \ref{eastheorem} and the examples in Section \ref{examplessection}.

\subsection{Acknowledgements}  The author is indebted to Vitaly Bergelson for encouragement and advice, and to Alexander Leibman for helpful comments on an early version of this article.

The author must also thank Michael Bj{\"o}rklund and Alexander Fish, whose discussions inspired the present work.

\section{Preliminaries from ergodic theory}\label{preliminaries}

\subsection{Measure preserving systems}

We make use of the theory of measure preserving systems, as presented in \cite{Gl} and \cite{W}.

For our purposes, a \textit{measure preserving system} $(X,\mathscr X,\mu, T)$ is a probability space $(X,\mathscr X,\mu)$, where $X$ is a set, $\mathscr X$ is a $\sigma$-algebra of subsets of $X,$ and $\mu$ is a probability measure on $\mathscr X,$ together with a transformation $T:X\to X$ such that $T^{-1}A\in \mathscr X$ whenever $A\in \mathscr X,$ and $\mu(T^{-1}A)=\mu(A)$ for all $A\in \mathscr X.$  We will always assume that $T$ is invertible.

If $\mathbf X=(X,\mathscr X,\mu,T)$ and $\mathbf Y=(Y,\mathscr Y,\nu,S)$ are measure preserving systems, we say that $\mathbf Y$ is a \textit{factor} of $\mathbf X$ if there exists a function $\pi:X\to Y$ satisfying $\pi^{-1}(\mathscr Y)\subset \mathscr X,$ $\mu(\pi^{-1}(D))=\nu(D)$ for all $D\in \mathscr Y,$ and $\pi(Tx)=S\pi(x)$ for $\mu$-almost every $x\in X.$  The map $\pi$ is called a \textit{factor map}.  If $\pi$ is one-to-one on a set of full measure, then we say that $\mathbf X$ and $\mathbf Y$ are isomorphic, and $\pi$ is an isomorphism.

Given a measure preserving system $(X,\mathscr X,\mu, T)$ with a factor $(Y,\mathscr Y,\nu, S)$ and factor map $\pi$, it is useful to consider the orthogonal projection $P_\mathbf Y:L^2(\mu)\to L^2(\mu)$ onto the space spanned by functions of the form $f\circ \pi,$ where $f\in L^2(\nu).$  This map can be identified with the conditional expectation $f\mapsto \mathbb E(f|\pi^{-1}(\mathscr Y)),$ so $P_\mathbf Y$ maps nonnegative functions to nonnegative functions.  See \cite{Fbook}, Chapter 5, for details and proofs.

A measure preserving system $(X,\mathscr X,\mu ,T)$ is called \textit{ergodic} if for all $A\in \mathscr X,$ $\mu(A \triangle T^{-1}A)=0$ implies that $\mu(A)=0$ or $\mu(A)=1.$  Equivalently, the system is ergodic if for all $f\in L^2(\mu), f\circ T=f$ implies that $f$ is $\mu$-almost everywhere equal to some constant.

\subsection{The ergodic theorem}\label{ergodicsection}

We will need the mean ergodic theorem and some of its consequences - see \cite{Gl} or \cite{RW} for proofs.

In this section, given a measure preserving system, $(X,\mathscr X,\mu, T),$ we let $P_I:L^2(\mu)\to L^2(\mu)$ denote the orthogonal projection onto the closed space of $T$-invariant functions.

\begin{theorem}\label{meanerg}  Let $\mathbf X=(X,\mathscr X,\mu,T)$ be a measure preserving system, and let $I\subset L^2(\mu)$ be the closed subspace of $T$-invariant functions.  Then for all $f\in L^2(\mu),$
\begin{align}\label{con}
\lim_{N-M\to \infty} \frac{1}{N-M}\sum_{n=M}^{N-1} f\circ T^n= P_I f
\end{align}
in the norm topology of $L^2(\mu).$
\end{theorem}
Note that if $\mathbf X$ is ergodic, then $P_I f=\int f\, d\mu.$

\subsubsection{Averaging sequences.}  We will need to consider averages more general than those in (\ref{con}). The following is shown in \cite{BE}, and in \cite{BKQW}, with slightly different terminology.

\begin{theorem}    Let $(\nu_j)_{j\in \mathbb N}$ be a sequence of probability measures on $\mathbb Z.$  The following conditions are equivalent.

\begin{enumerate}
\item[(i)]  For all $\theta\in (0,2\pi),$ $\lim_{j\to \infty} \int \exp(in\theta)\, d\nu_j(n)=0.$
\item[(ii)]  For every measure preserving system $(X,\mathscr X,\mu,T)$ and all $f\in L^2(\mu),$
$$
\lim_{j\to \infty} \int f\circ T^n d\nu_j(n) = P_If,
$$ in the norm topology of $L^2(\mu).$
\end{enumerate}
\end{theorem}

As a consequence, given a function $f\in L^2(\mu)$ and a sequence $(\nu_j)_{j\in \mathbb N}$ satisfying (i) or (ii) above, one may pass to a subsequence $(\nu_j')_{j\in \mathbb N}$ to obtain pointwise $\mu$-almost everywhere convergence of the averages  $\int f\circ T^n \, d\nu_j'(n).$

We call a sequence $(\nu_j)_{j\in \mathbb N}$ satisfying (i) or (ii) above an \textit{equidistributed averaging sequence}.

The mean ergodic theorem says that the sequence of measures defined by $\nu_j:=$ normalized counting measure on $\{1,\dots, j\}$ is an equidistributed averaging sequence.  Many examples of sparsely supported equidistributed averaging sequences are given by Theorems 3.2 and 3.3 of \cite{BKQW};  here are two special cases.  We use $\lfloor x\rfloor$ to denote the greatest integer less than or equal to $x.$
\begin{itemize}
 \item[$\bullet$]  $\nu_j=$ normalized counting measure on $\{\lfloor n^\alpha\rfloor: 1\leq n \leq j\},$ for $0<\alpha \notin \mathbb Z.$
 \item[$\bullet$]  $\nu_j=$ normalized counting measure on $\{\lfloor \sqrt{2}n^4-\pi n^2 \rfloor: 1\leq n \leq j\}.$
\end{itemize}

As we shall see in Section \ref{examplessection}, the set $\{n^2: n\in \mathbb N\}$ does not support an equidistributed averaging sequence, nor does the set of primes.

\subsection{Kronecker systems and the Kronecker factor}

For our purposes, a \textit{group rotation} is a measure preserving system $\mathbf Z=(Z,\mathscr Z, m_Z, R_\alpha),$ where $Z$ is a compact abelian group, $\mathscr Z$ is the Borel $\sigma$-algebra of $Z, m_Z$ is Haar measure, $\alpha\in Z,$ and $R_\alpha$ is defined by $R_\alpha(z)=z+\alpha$ for $z\in Z.$  We do not assume that $Z$ is metric.

If $\{n\alpha: n\in \mathbb Z\}$ is dense in $Z$ then $\mathbf Z$ is ergodic, and we say that $\mathbf Z$ is a \textit{Kronecker system}.  It is well known that such systems are minimal (as topological systems), and hence that for all nonempty open $U\subset Z$ and $z\in Z,$ the set of entry times $\{n:z+n\alpha\in U\}$ is syndetic.

 If $\mathbf X$ is an ergodic measure preserving system, there is a factor $\mathbf Y=(Y,\mathscr Y,\nu,S)$ with factor map $\pi$ having the property that the eigenfunctions of $T$ (that is, those $f\in L^2(\mu)$ satisfying $f\circ T=\lambda f$ $\mu$-almost everywhere, for some $\lambda\in \mathbb C$) are measurable with respect to $\pi^{-1}(\mathscr Y),$ and $\mathbf Y$ is the smallest such factor in the sense that $\pi^{-1}(\mathscr Y)$ is generated by the eigenfunctions of $T.$  This factor is unique up to isomorphism, and is called the \textit{Kronecker factor} of $\mathbf X.$  Since $\mathbf Y$ is ergodic and $L^2(\nu)$ is spanned by the eigenfunctions of $S,$  the Halmos-von Neumann theorem (see \cite{KFS} or \cite{Gl}) says that $\mathbf Y$ is isomorphic to a compact group rotation $(Z,\mathscr Z,m,R_\alpha).$  Given a system $\mathbf X,$ we will denote its Kronecker factor by $\mathbf Z,$ and we will assume that $\mathbf Z$ is actually a compact group rotation, not merely that it is isomorphic to such a system.  Also, we will abuse notation and write $\mathbb E(f|\mathscr Z)$ for $\mathbb E(f|\pi^{-1}(\mathscr Z)).$

The following classical result describes the set of $T\times T$-invariant functions, given an ergodic system $(X,\mathscr X,\mu,T).$  See \cite{Gl}, Chapter 9, for a more general result and proof.

\begin{lemma}\label{invariant}
If $(X,\mathscr X,\mu,T)$ is an ergodic system and $(Z,\mathscr Z,m, R_\alpha)$ is its Kronecker factor, with factor map $\pi,$ then the space of $T\times T$-invariant functions in $L^2(\mu\times \mu)$ is contained in the closed space spanned by functions of the form $(x,y)\mapsto f(\pi(x))g(\pi(y)),$ with $f,g\in L^2(m).$
\end{lemma}

\subsubsection{Ergodic averages for Kronecker systems and limits of measures.}  Ergodic averages on Kronecker systems are particularly well behaved.

\begin{lemma}\label{uniform}  Let $(Z,\mathscr Z,m, R_\alpha)$ be a Kronecker system, let $f:Z\to \mathbb C$ be continuous, and let $(\nu_j)_{j\in \mathbb N}$ be an equidistributed averaging sequence.  Then the averages $\int f(z+n\alpha) d\nu_j(n)$ converge uniformly to $\int f\, dm$ as $j\to \infty.$
\end{lemma}

\textit{Proof.}  It suffices to show that the conclusion holds when $f$ is a character of $Z,$ since the characters of $Z$ span a uniformly dense subspace of $C(Z).$  This special case follows from the definition of equidistributed averaging sequence, since for a character $\chi:Z\to \mathbb C,$ we have
$$
\int \chi(z+n\alpha)\, d\nu_j(n)= \chi(z) \int \chi(\alpha)^n\, d\nu_j(n)
$$
for all $z\in Z.$  We know that the limit is $\int \chi\, dm,$ since the limit in $L^2(\mu)$ is $\int \chi\, dm.$  \hfill $\square$

\medskip

We will also consider averages
\begin{align}\label{convos}
\int f\circ R_\alpha^n\, d\eta_j(n)
\end{align} where $f\in L^\infty(Z)$ and $\eta_j\leq \nu_j$ in the sense that $\eta_j(\{n\})\leq \nu_j(\{n\})$ for each $n\in \mathbb Z.$ The averages (\ref{convos}) are simply linear combinations of the functions $f\circ R_\alpha^n,$ although we may interpret them as convolutions of measures.

The next lemma describes the weak limits of the averages in (\ref{convos}).  If $Z$ is a compact metric abelian group with Haar measure $m$ and $f,g\in L^\infty(m),$ we consider the convolution $f*g$ defined by $f*g(t):=\int f(z)g(t-z)\, dm(z).$

\begin{lemma}\label{weak*limits}  Let $(Z,\mathscr Z, m,R_\alpha)$ be a Kronecker system (with $Z$ metrizable), and let $(\nu_j)_{j\in \mathbb N}$ be an equidistributed averaging sequence.  Suppose that $\eta_j(\{n\})\leq \nu_j(\{n\})$ for all $n\in \mathbb Z,$ while $\lim_{j\to\infty} \eta_j(\mathbb Z)=c>0.$  If $f\in L^\infty(m)$ and $\varphi\in L^2(m)$ is a weak limit of $\int f\circ R_\alpha^{-n} \, d\eta_j(n)$ in the sense that
$$
\lim_{j\to \infty} \int \int f\circ R_\alpha^{-n} \, d\eta_j(n)\cdot h\, dm= \int \varphi \cdot h\, dm
$$
for all $h\in L^2(m),$  then $\varphi=f*\psi$ for some measurable $\psi:Z\to [0,1]$ with $\int \psi\, dm=c.$
\end{lemma}

\textit{Proof.}  For $j\in \mathbb N,$ let $\eta_j^*$ be the measure on $Z$ given by $\int h\, d\eta_j^*= \int h(n\alpha) d\eta_j(n)$ for continuous $h:Z\to \mathbb R.$  Passing to a subsequence, we may assume that the $\text{weak}^*$ limit of the $\eta_j^*$ exists; call the limit $\eta^*.$  We claim that $\eta^*$ is absolutely continuous with respect to Haar measure and its Radon-Nikodym derivative $\psi$ is bounded above by $1.$  This follows from the fact that $\int h\, d\eta_j^*\leq \int h \, d\nu_j^*$ for every continuous positive function $h\in C(Z),$ while $\lim_{j\to \infty} \int h \, d\nu_j^*= \int h\, dm,$ by Lemma \ref{uniform};  hence $\eta^*(K)\leq m(K)$ for every compact $K\subset Z.$  For this $\psi:Z\to [0,1]$  and every $h\in C(Z),$ we now have
\begin{align}\label{RND}
\lim_{j\to \infty} \int h\, d\eta_j^*=\int h\, d\eta^*=\int h\cdot \psi\, dm.
\end{align}
Furthermore, $\int \psi \, dm=c,$ since $\lim_{j\to \infty} \eta_j(\mathbb Z)=c.$

Note that it suffices to establish the lemma when $f$ is a character of $Z,$ since the characters span a dense subset of $L^2(m).$  With this assumption, we have, for all $z\in Z,$
\begin{align*}
\lim_{j\to\infty} \int f\circ R_\alpha^{-n}(z)\, d\eta_j(n)&= \lim_{j\to\infty} f(z)\int f(-w)\, d\eta_j^*(w) \\
&=f(z)\int f(-w)\, d\eta^*(w)\\
&=f(z) \int f(-w)\cdot \psi(w)\, dm(w),
\end{align*}
the last equality being an instance of (\ref{RND}).  Since 
$$
f(z) \int f(-w)\cdot \psi(w)\, dm(w)=\int f(z-w)\cdot \psi(w)\, dm(w)=f*\psi(z),
$$ this completes the proof.  \hfill $\square$

\subsection{The Steinhaus lemma}

We formulate a convenient version of the Steinhaus lemma for compact abelian groups.

\begin{lemma}\label{Steinhaus+} Let $Z$ be a compact abelian group with Haar measure $m,$ and let $f,g:Z\to [0,1]$ be measurable functions with $\int f\, dm>0, \int g\, dm>0.$  Then $f*g$  is continuous, and its support has measure at least $\max\{\int f\, dm,\int g\, dm\}.$
\end{lemma}
Here ``the support of $h$" means $\{x:h(x)>0\}.$  In particular, if $f=1_C,g=1_D$ for measurable sets $C,D\subset Z,$ then $f*g$ is supported on $C+D,$ so $C+D$ contains an open set with measure at least $\max\{m(C),m(D)\}.$

\medskip

To prove Lemma \ref{Steinhaus+}, expand $f$ and $g$ as Fourier series $f=\sum_{\chi\in \widehat{Z}} \hat{f}(\chi) \chi,$ $g= \sum_{\chi \in \hat Z} \hat{g}(\chi) \chi,$ and note that $\hat{f},\hat{g}\in L^2(\widehat{Z}).$  Then $\hat{f}\cdot \hat{g}\in L^1(\widehat{Z}),$ so $f*g=\sum_{\chi \in \widehat{Z}} \hat{f}(\chi)\hat{g}(\chi)\chi$ is a uniform limit of continuous functions, and so is continuous.  To estimate the support of $f*g,$ note that $\int f*g \, dm= \int f\, dm \int g\, dm$ by Fubini's theorem, while $\sup_{t\in Z} f*g(t)\leq \min\{\int f\, dm, \int g\, dm\},$ so
\begin{align*}
m\{t: f*g(t)>0\}\cdot \min\{\int f\, dm, \int g\, dm\} \geq \int f\, dm \int g\, dm
\end{align*}
The last inequality implies $m\{t: f*g(t)>0\}\geq \max\{\int f\, dm, \int g\, dm\}.$

\medskip

\begin{remark}While the preceding proof is standard, it exhibits a theme in common with the work on sumsets in finite groups mentioned in Section \ref{Introduction}: in the setting of $\mathbb Z/N\mathbb Z,$ one considers the characteristic functions $f,g$ of two sets $A,B\subset \mathbb Z/N\mathbb Z,$ and then uses bounds on the $L^2$ norm of $\hat f$ and $\hat g$ to obtain a bound on the $L^1$ norm of $\widehat{f*g}.$  With much effort, this bound is exploited to reveal the structure of the support of $f*g,$ and hence the structure of $A+B.$
\end{remark}

We now consider a partial converse to the Steinhaus lemma; we need it to construct the examples in Section \ref{examplessection}.

\begin{lemma}\label{Steinhausconverse}  Let $Z$ be a separable compact abelian group with Haar measure $m,$ and let $E\subset Z$ be compact with $m(E)=0.$  For all $\varepsilon>0,$ there exists a compact $K\subset Z$ with $m(K)>1-\varepsilon$ such that $E+K$ has empty interior.
\end{lemma}

\textit{Proof.}  Fix $\varepsilon>0.$  Let $\{V_n:n\in \mathbb N\}$ be a collection of open sets whose union is dense in $Z$ such that $m(V_n-E)<\varepsilon 2^{-n}$ for each $n.$  Let $K=\bigcap_n Z\setminus (V_n-E).$  Then $m(K)\geq 1-\sum_{n}\varepsilon 2^{-n}= 1-\varepsilon,$ and $K$ is compact.  Furthermore $(E+K)\cap V_n=\emptyset$ for all $n,$ so the complement of $E+K$ is dense, hence $E+K$ has empty interior.  \hfill $\square$

\subsection{Bohr sets in $\mathbb Z$}\label{bs}

The Bohr topology on $\mathbb Z$ is the topology generated by the functions $n\mapsto \exp(i\lambda n), \lambda \in \mathbb R.$  A basis for this topology consists of sets of the form $\{n:\operatorname{Re} p(n)>0\},$ where $p$ is a trigonometric polynomial given by $p(n)=\sum_{\lambda \in F} c_\lambda \exp(i\lambda n)$ for some $c_\lambda\in \mathbb C$ and finite $F\subset \mathbb R.$  A set $S\subset \mathbb Z$ is called a \textit{Bohr set} if it contains one of these nonempty basis sets.

Equivalently, we call $B\subset \mathbb Z$ a Bohr set if it contains a set of entry times to an open set in a compact metric Kronecker system.  That is, $B$ is a Bohr set if there exists an ergodic group rotation system $(Z,\mathscr Z, m, R_\alpha)$ and an open $U\subset Z$ such that $B$ contains $\{n:n\alpha\in U\}.$  The equivalence of the two definitions of ``Bohr set" follows from Pontryiagin duality.

Following \cite{BFW}, we call a set \textit{piecewise Bohr} if it contains the intersection of a Bohr set and a thick set.   Such a set is piecewise syndetic, and an example given in \cite{BFW} shows that there are syndetic sets that are not piecewise Bohr.  Another such example is the set 
$$
S:=\{n: n^2\sqrt{2} \mod 1 \in (0,1/2)\};
$$ Weyl's theorem on equidistribution implies that $S$ is syndetic, and the same theorem implies that for all sets of form $R=\{n: \operatorname{Re} p(n)>~0\}$ where $p:\mathbb Z\to \mathbb C$ is a trigonometric polynomial, 
$$
d^*(S\cap R)=\frac{1}{2}\lim_{N-M\to \infty} \frac{|R\cap [M,N]|}{N-M+1}.
$$  Hence $S$ cannot contain the intersection of such an $R$ with a thick set, and so cannot be piecewise Bohr.

One may view piecewise Bohr sets as having more structure than arbitrary sets of positive density.  For instance, it is shown in \cite{BHK} that for every $m>0,$ there are sets $B\subset \mathbb Z$ with $d^*(B)>0,$ while $d^*\bigl(\bigcap_{l=0}^4 B-ld\bigr)< d^*(B)^m/2$ for all $d\neq 0.$  In other words, $B$ contains much less than the expected density of $5$-term arithmetic progressions of a given common difference $d,$ for every $d\neq 0.$  In contrast, a given piecewise Bohr set will have, for many $d,$ more than the expected density of $k$-term arithmetic progressions with difference $d.$  To be more precise, we state the following lemma.

\begin{lemma}\label{excess}  Let $(Z,\mathscr Z,m,R_\alpha)$ be an ergodic Kronecker system, and let $U\subset Z$ be open.  Let $B_0=\{n:n\alpha\in U\},$ and let $B$ be the intersection of $B_0$ with a thick set.  Then for all $k\in \mathbb N$ and all $\varepsilon>0,$ there exists $d\in \mathbb N$ such that
\begin{align}\label{bigint}
d^*\bigl(\bigcap_{l=0}^k B-ld\bigr)>m(U)-\varepsilon.
\end{align}
In fact, there is a Bohr set of $d$ satisfying \rm{(\ref{bigint})}.
\end{lemma}

\textit{Proof.}  Fix $k\in \mathbb N,$ and write $W_d$ for $\bigcap_{l=0}^k U-ld\alpha.$ Choose $d\in \mathbb Z,$ with $m(W_d)>m(U)-\varepsilon$ (note that there is a Bohr set of such $d$). Let $f:Z\to [0,1]$ be a continuous function supported on $W_d$ with $\int f\, dm> m(U)-\varepsilon.$  Let $(I_r)_{r\in \mathbb N}$ be a sequence of intervals with $|I_r|\to \infty$ and $I_r\cap B_0\subset B$ for all $r.$  Then $\bigcap_{l=0}^k B-ld$ contains $\bigcap_{l=0}^k (I_r\cap B_0)-ld$ for each $r.$  We will show that $ |I_r\cap \bigcap_{l=0}^k (I_r\cap B_0)-ld|/|I_r|>m(U)-\varepsilon$ for sufficiently large $r.$

By Lemma \ref{uniform}, we have
\begin{align}\label{detect}
\lim_{r\to \infty} \frac{1}{|I_r|}\sum_{n\in I_r} f(n\alpha)=\int f\, dm,
\end{align}
and the integral is at least $m(U)-\varepsilon.$  If $f(n\alpha)>0,$ then $n\alpha\in U-ld\alpha$ for $0\leq l\leq k,$ meaning $n\in B_0-ld$ for each $l.$  Since $f$ is supported on $W_d,$ (\ref{detect}) implies that $\liminf_{r\to \infty}  |I_r\cap  \bigcap_{l=0}^k B_0-ld|/|I_r|\geq m(U)-\varepsilon,$ and the conclusion follows. \hfill $\square$

\medskip

As observed in \cite{BFW} and \cite{BBF},  if a Bohr set is ``cut into long segments, shifted, and reassembled," the resulting set is piecewise Bohr.   The next lemma makes this statement precise.

\begin{lemma}\label{shifts} Let $B\subset \mathbb Z$ be a Bohr set, let $(I_j)_{j\in \mathbb N}$ be a sequence of intervals with lengths tending to infinity, and let $(r_j)_{j\in \mathbb N}$ be a sequence of integers.  Then $\bigcup_j (I_j\cap B)+r_j$ is piecewise Bohr.
\end{lemma}

\textit{Proof.}  Let $(Z,\mathscr Z,m,R_\alpha)$ be a metric Kronecker system,  with $U\subset Z$ open such that $\{n: n\alpha\in U\}$ is contained in $B.$  Since $Z$ is compact metric, there is a subsequence $(r_j'\alpha)_{j\in \mathbb N}$ of $(r_j\alpha)_{j\in \mathbb N}$ that converges to a point $z_0\in Z.$ Choosing $J$ sufficiently large, $\bigcap_{j>J} U+r_j'\alpha$ contains an open set $V.$

Now for sufficiently large $j,$ $(B\cap I_j')+r_j'$ contains $\{n\in I_j'+r_j':n\alpha \in U+r_j'\},$ so $(B\cap   I_j')+r_j'$ contains $\{n\in I_j'+r_j': n\alpha\in V\}.$  Thus $\bigcup_j (I_j\cap B)+r_j$ contains the intersection of the thick set $\bigcup_{j>J} I_j'+r_j'$ with the Bohr set $\{n:n\alpha \in V\}.$  \hfill $\square$

\medskip

The proof of Lemma \ref{shifts} gives a bit more information; we could choose $r_j'$ above so that $m(\bigcap_{j} U+r_j'\alpha)>m(U)-\varepsilon.$  This leads to the following refinement.

\begin{corollary}\label{shifts+}  Let $(Z,\mathscr Z,m, R_\alpha)$ be a metric Kronecker system with $U\subset Z$ open. Let $B=\{n\in \mathbb Z: n\alpha \in U\}$ , let $(I_j)_{j\in \mathbb N}$ be a sequence of intervals with lengths tending to infinity, and let $(r_j)_{j\in \mathbb N}$ be a sequence of integers.  Then $\bigcup_j (I_j\cap B)+r_j$ contains the intersection of a thick set with some $B':=\{n\in \mathbb Z: n\alpha\in V\},$ where $V$ is open and $m(V)> m(U)-\varepsilon.$
\end{corollary}

\subsection{The Bohr compactification}\label{bc}  Although not strictly necessary for our proofs and examples, the Bohr compactification of $\mathbb Z$ provides a useful perspective for some of the questions asked in Section \ref{questions}.

To form the Bohr compactification $b\mathbb Z$ of $\mathbb Z,$ give $\mathbb R/\mathbb Z$ the discrete topology, and let $b\mathbb Z$ be the dual of that discrete group.  Then $b\mathbb Z$ is compact, and $\mathbb Z$ embeds densely therein by $n\mapsto e_n,$ where $e_n(t)=\exp(2\pi i n t).$  Under this embedding, the characters of $\mathbb Z$ extend continuously to characters of $b\mathbb Z$, and the Bohr topology on $\mathbb Z$ is the subspace topology on $\mathbb Z$ induced by the topology on $b\mathbb Z.$  See \cite{R} for details.

In the sequel, we will consider $\mathbb Z$ as a subset of $b\mathbb Z,$ and we can speak of a sequence of measures on $\mathbb Z$ converging in the $\text{weak}^*$ topology of $b\mathbb Z.$  In particular, a sequence $(\nu_j)_{j\in \mathbb N}$ of probability measures on $\mathbb Z$ converges to Haar measure on $b\mathbb Z$ if and only if it is an equidistributed averaging sequence.

\subsection{A weak correspondence principle}  We will use  a weak version of Furstenberg's correspondence principle from \cite{F77}.  To state our version, we fix notation for the shift space.  Let $\Omega=\{0,1\}^\mathbb Z$ with the product topology, and define the shift $\sigma:\Omega\to \Omega$ by $(\sigma x)(n)=x(n+1)$ for $x\in \Omega.$  Then $\Omega$ is a compact metric space and $\sigma$ is a surjective homeomorphism.

\begin{proposition}\label{correspondence}  Suppose $B\subset \mathbb Z$ with $d^*(B)>0.$  Let $X=\overline{\{\sigma^n 1_B:n\in \mathbb Z\}},$ the orbit closure of $1_B$ in the shift space $(\{0,1\}^\mathbb Z,\sigma)$ and let $O{}$ be the open set $\{x\in~X: x(0)=1\}.$  Then there is a $\sigma$-invariant probability measure $\mu$ on $X$ with $\mu(O{})\geq d^*(B).$  Furthermore, we can pick $\mu$ so that $(X,\mathscr X,\mu,\sigma)$ is ergodic.
\end{proposition}

\textit{Proof.}  For $B\subset \mathbb Z$, let $x=1_B,$ and let $I_k$ be a sequence of intervals with $|I_k|\to \infty$ and $\lim_{k\to\infty} \frac{|B\cap I_k|}{|I_k|}=d^*(B).$  Let $\mu_k=\frac{1}{|I_k|}\sum_{n\in I_k}^N \delta_{\sigma^nx},$ where $\delta_{\sigma^nx}$ is the unit point mass at $\sigma^n x.$  Then every $\text{weak}^*$ limit $\mu$ of the sequence $(\mu_k)_{k\in \mathbb N}$ is $\sigma$-invariant and satisfies $\mu(O{})\geq d^*(B).$  We can find an ergodic $\mu$ with these properties by applying the ergodic decomposition theorem (see \cite{Gl}).  \hfill $\square$

\medskip

Proposition \ref{correspondence} also follows from the proof of Proposition 3.1 of \cite{BHK}.

\section{Proof of Theorems \ref{eastheorem} and \ref{APs}}\label{ergodicschemes}

\subsection{Refinement and ergodic theoretic analogue}

We will deduce Theorems \ref{eastheorem} and \ref{APs} from the next theorem, which is a refinement of Theorem \ref{eastheorem}.

\begin{theorem}\label{eastheorem'}  Let $(\nu_j)_{j\in \mathbb N}$ be an equidistributed averaging sequence, and let $A,B\subset \mathbb Z$ with $d^*(B)>0.$  Then the following implications hold.
\begin{enumerate} \item[1.]  If $d_\nu(A)>0,$ then $A+B$ is piecewise Bohr.  In fact, there exists a Kronecker system $(Z,\mathscr Z,m,R_\alpha)$ such that for all $\varepsilon>0,$ there is an open set $U\subset Z$ with $m(U)>\max\{d_\nu(A),d^*(B)\}-\varepsilon$ and a thick set $S\subset \mathbb Z$ such that $A+B$ contains $S\cap \{n\in \mathbb Z: n\alpha\in U\}.$
\item[2.]  If $d_\nu(A)=1,$ then $A+B$ is thick.
\end{enumerate}
\end{theorem}

We will deduce Theorem \ref{eastheorem'} from the following proposition about measure preserving systems.

\begin{proposition}\label{Kroneckeropen}  Let $(X,\mathscr X,\mu, T)$ be an ergodic measure preserving system with Kronecker factor $(Z,\mathscr Z,m, R_\alpha)$ and factor map $\pi:X\to Z.$ Let $D\in \mathscr X$ with $\mu(D)>0,$ and let $(\nu_j)_{j\in \mathbb N}$ be an equidistributed averaging sequence.  If $A\subset \mathbb Z$ with $d_\nu(A)>0,$ then $\bigcup_{a\in A} T^{a}D$ contains, up to $\mu$-measure $0,$ a set of the form $\pi^{-1}(U),$ where $U\subset Z$ is open and $m(U)\geq \max\{d_\nu(A),\mu(D)\}.$
\end{proposition}
Postponing the proof of Proposition \ref{Kroneckeropen} to Subsection \ref{postponed}, we proceed with the proof of Theorem \ref{eastheorem'}.

\medskip

\textit{Proof of Theorem \ref{eastheorem'}.}  Fix an equidistributed averaging sequence $(\nu_j)_{j\in \mathbb N},$ and sets $A,B\subset \mathbb Z$ with $d_\nu(A)>0$ and $d^*(B)>0.$  Let $T$ be the shift on $\{0,1\}^{\mathbb Z},$ let $X$ be the orbit closure of $1_B$ in $(\{0,1\}^{\mathbb Z},T),$ and let $O=\{x\in X:x(0)=1\}.$   As Proposition \ref{correspondence} allows, let $\mu$ be a $T$-invariant probability measure on $X$ so that $(X,\mathscr X,\mu, T)$ is ergodic and $\mu(O{})\geq d^*(B).$

\medskip

Define $V:=\bigcup_{a\in A} T^{a}O{}.$

\begin{lemma}\label{ergodictosumset}
With $A,B,X,$ and $V$ as above, let $x\in X.$ For all finite sets $F\subset \{n:T^nx\in V\},$ $A+B$ contains a translate of $F.$
\end{lemma}

\textit{Proof.}  Since $x$ is in the $T$-orbit closure of $1_B,$ we can write $x=1_E,$ where $E\subset \mathbb Z$ is an \textit{increasing} union of the form
$$
\bigcup_{k=1}^\infty ([n_k-k, n_k+k]\cap B)-n_k,
$$ for some sequence of integers $n_k.$  This means that $A+E$ is an increasing union $\bigcup_{k=1}^\infty A+([n_k-k, n_k+k]\cap B)-n_k.$  In particular, for a given finite interval $I$,
$$
(A+E)\cap I= (A+([n_k-k,n_k+k]\cap B)-n_k)\cap I \text{\  for some $k$,}
$$
hence $(A+B-n_k)\cap I$ contains $(A+E)\cap I.$

It follows from the definitions of $E$ and $V$ that $A+E$ is simply $\{n:T^nx\in V\},$ so the Lemma is proved. $\square$

\medskip

We now complete the proof of Theorem \ref{eastheorem'}, Part 1. Let $(Z,\mathscr Z,m,R_\alpha)$ be the Kronecker factor of $(X,\mathscr X,\mu, T),$ and let $\pi:X\to Z$ be the factor map.  Fix $\varepsilon>0.$  By Proposition \ref{Kroneckeropen}, $V$ ($=\bigcup_{a\in A} T^aO{}$) contains, up to $\mu$-measure $0,$ a set of the form $\pi^{-1}(U),$ where $U\subset Z$ is open and $m(U)>\max\{d_\nu(A),d^*(B)\}-\varepsilon.$  Let $x\in X$ such that
\begin{enumerate}
\item[(i)]  The equation $\pi(T^nx)=\pi(x)+n\alpha$ holds for all $n\in \mathbb Z,$ and
\item[(ii)]  $\pi(T^nx)\in U$ implies that $T^nx\in V.$
\end{enumerate}
Since the sets of $x$ satisfying each of (i) and (ii) separately have full measure, there is an $x$ satisfying both conditions.  To see that the set of $x$ satisfying (ii) has full measure, note that its complement is contained in $\bigcup_{n\in \mathbb Z} T^{-n} (\pi^{-1}(U)\setminus V).$

By condition (i), the set $\{n:\pi(T^nx)\in U\}$ is a Bohr set, and by condition (ii), this set is contained in $\{n:T^nx\in V\}.$  Since $m(U)\geq \max\{d_\nu(A),d^*(B)\}-\varepsilon,$ Lemma \ref{ergodictosumset} and Corollary \ref{shifts+} now imply Part 1 of Theorem \ref{eastheorem'}.

Proof of Theorem \ref{eastheorem'}, Part 2.  Without loss of generality we can assume that every $\nu_j$ is supported on $A.$  Let $f=1_{O{}},$ and let $C\in \mathscr X$ with $\mu(C)>0.$  Since $(\nu_j)_{j\in \mathbb N}$ is an equidistributed averaging sequence, the averages
$$
\int \int f\circ T^{-n}\cdot 1_C \, d\mu\, d\nu_j(n)
$$
converge to $\int f\, d\mu \int 1_C\, d\mu.$  In particular, there exists $n\in A$ such that $\mu(E\cap T^n O{})>0.$  Since this is true for every $C$ of positive measure, it follows that $V=\bigcup_{n\in A} T^nO{}$ has full measure.  Hence, there exists $x\in X$ such that $T^n x\in V$ for all $n\in \mathbb Z.$  Lemma \ref{ergodictosumset} now implies that $A+B$ contains a shift of every finite subset of $\mathbb Z,$ and in particular that $A+B$ contains intervals of every finite length.  \hfill $\square$

\medskip

\textit{Proof of Theorem \ref{APs}.}  Theorem \ref{APs} follows directly from Part 1 of Theorem \ref{eastheorem'}, together with Lemma \ref{excess}.  \hfill $\square$

\begin{remark} The deduction of Theorem \ref{eastheorem'} from Proposition \ref{Kroneckeropen} is similar to the proof of \cite{Fbook}, Theorem 3.20, a result of R.~Ellis which says that every set $B$ with $d^*(B)>0$ contains translates of every finite subset of some set $B'$ having density $d^*(B)$.  That is, $d(B'):=\lim_{N\to \infty} \frac{|B\cap [1,N]|}{N}=d^*(B),$ and for every finite $F\subset B',$ there exists $c$ with $F+c\subset B.$
\end{remark}

\subsection{Proof of Proposition \ref{Kroneckeropen}.}\label{postponed}

Let us recall Proposition \ref{Kroneckeropen} and describe the idea of the proof.

\begin{proposition 3.4} Let $(X,\mathscr X,\mu, T)$ be an ergodic measure preserving system with Kronecker factor $(Z,\mathscr Z,m, R_\alpha)$ and factor map $\pi:X\to Z.$ Let $D\in \mathscr X$ with $\mu(D)>0,$ and let $(\nu_j)_{j\in \mathbb N}$ be an equidistributed averaging sequence.  If $A\subset \mathbb Z$ with $d_\nu(A)>0,$ then $\bigcup_{a\in A} T^{a}D$ contains, up to $\mu$-measure $0,$ a set of the form $\pi^{-1}(U),$ where $U\subset Z$ is open and $m(U)\geq \max\{d_\nu(A),\mu(D)\}.$
\end{proposition 3.4}

To prove the Proposition, we will bound $1_{\bigcup_{a\in A}T^a D}$ from below by averages $g_j:=\int 1_D\circ T^{-n} \, d\eta_j(n),$ where $\eta_j(E):=\nu_j(A\cap E);$ one easily verifies $g_j\leq 1_{\bigcup_{a\in A} T^a D}.$  We then pass to a subsequence to obtain a weak limit $g:=\lim_{j\to \infty} g_j.$  With the aid of the next lemma, we find that $g$ is equal to $\lim_{j\to \infty} \int \mathbb E(1_D|\mathscr Z)\circ T^{-n} \, d\eta_j(n).$   Thinking of $\mathbb E(1_D|\mathscr Z)$ as a function $f:Z\to [0,1],$ we use Lemma \ref{weak*limits} to describe $g.$

\medskip

The next lemma is standard in multiple recurrence arguments; cf.~\cite{Fbook}, Lemma 4.15.

\begin{lemma}\label{wm}
Suppose that $(X,\mathscr X,\mu, T)$  is an ergodic measure preserving system with Kronecker factor $(Z,\mathscr Z, m, R_\alpha).$  Suppose that $f\in L^2(\mu)$ with $\mathbb E(f|\mathscr Z)=0,$ and let $(\nu_j)_{j\in \mathbb N}$ be an equidistributed averaging sequence.  Then for all $g\in L^2(\mu)$ and all $\varepsilon>0,$
$$
d_\nu\Bigl\{n:\bigl| \int f\circ T^{-n}\cdot g \, d\mu\bigr |>\varepsilon\Bigr\}=0.
$$
\end{lemma}

\textit{Proof.}  Write $P_Zf$ for $\mathbb E(f|\mathscr Z).$  The conclusion is equivalent to the assertion that
$$
\lim_{j\to \infty} \int \bigl| \int f\circ T^{-n} \cdot g\, d\mu \bigr|^2 \, d\nu_j(n)=0.
$$

Writing $|\int f\circ T^{-n} \cdot g\, d\mu|^2=\int f\otimes \bar f\circ  (T\times T)^{-n} \cdot g\otimes \bar g\, d\mu\times \mu,$ we average with respect to $\nu_j$ to find
\begin{align}\label{wmavg}
\lim_{j\to \infty} \int \bigl|\int f\circ T^{-n} \cdot g\, d\mu\bigr|^2\, d\nu_j(n)= \int P_{T\times T} (f\otimes \bar f)\cdot g\otimes \bar g \, d\mu\times \mu,
\end{align}
where  $P_{T\times T} (f\otimes \bar f)$ is the projection of $f\otimes \bar f$ on the space of $T\times T$-invariant functions in $L^2(\mu\times \mu).$  By Lemma \ref{invariant}, the space of $T\times T$-invariant functions in $L^2(\mu\times \mu)$ is spanned by functions of the form $P_Z(h_1)\otimes P_Z(h_2).$  Since $P_Z(f)=0,$ the funtcion $f\otimes f$ is orthogonal to the space of $T\times T$-invariant functions.  Hence, the integral on the right-hand side of (\ref{wmavg}) is $0.$  \hfill $\square$

\medskip

\textit{Proof of Proposition \ref{Kroneckeropen}.}  Let $(X,\mathscr X,\mu,T),$ $(Z,\mathscr Z, m, R_\alpha)$ and $\pi:X\to Z$ be as in the hypotheses of the proposition.  Passing to a subsequence of $(\nu_j)_{j\in \mathbb N},$ we suppose that $\lim_{j\to\infty} \nu_j(A)$ exists and equals $d_\nu(A).$   Consider the measures $\eta_j$ on $\mathbb Z$ given by $\eta_j(E)=\nu_j(A\cap E.)$  Then for all $j, g_j:=\int 1_D\circ T^{-n} \, d\eta_j(n)$ is supported on $\bigcup_{a\in A} T^{a} D,$ as is any weak limit of the $g_j.$  Again passing to a subsequence, we may assume that $g:=\lim_{j\to \infty} g_j$ exists weakly, in the sense that $\lim_{j\to \infty} \int g_j \cdot h\, d\mu$ exists for all $h\in L^2(\mu).$ Thus Proposition \ref{Kroneckeropen} is a consequence of the following claim.

\begin{claim} The support of $g:=\lim_{j\to \infty} g_j$ contains a set of the form $\pi^{-1}(U),$ where $U\subset Z$ is open and $m(U)\geq \max\{d_\nu(A), \mu(D)\}.$
\end{claim}

Note: we use ``support of $g$" to mean $\{x:g(x)>0\}.$

To prove the claim write $1_D=f_1+f_0,$ where $f_1=\mathbb E(1_D|\mathscr Z), \mathbb E(f_0|\mathscr Z)=0.$  Then $g_j$ decomposes as $g_{j,1}+g_{j,0},$ where $g_{j,1}=\int f_1\circ T^{-n} d \eta_j(n),$\\ $g_{j,0}= \int f_0\circ T^{-n}\, d\eta_j(n).$

By Lemma  \ref{wm}, we have for all $h\in L^2(\mu)$ and all $\varepsilon>0,$
$$
d_\nu\Bigl\{n:\bigl|\int f_0\circ T^{-n}\cdot h \, d\mu\bigr|>\varepsilon\Bigr\}=0.
$$  From this and the fact that $d_\nu(A)>0$ we conclude that
\begin{align*}
d_\eta\Bigl\{n:\bigl|\int f_0\circ T^{-n}\cdot h \, d\mu\bigr|>\varepsilon\Bigr\}=0.
\end{align*}
  Hence for all $h\in L^2(\mu),$
$$
\lim_{j\to \infty} \int g_{j,0} \cdot h \, d\mu= \lim_{j\to \infty} \int \int f_0\circ T^{-n}\cdot h\, d\mu\, d\eta_j(n)=0.
$$
Thus $\lim_{j\to \infty} g_j=\lim_{j\to \infty} g_{j,1}.$  Since $f_1$ is $\mathscr Z$-measurable, write $f_1= \tilde f_1\circ \pi,$ where $\tilde f_1\in L^\infty(Z,m),$ and similarly write $g_{j,1}={\tilde g}_{j,1}\circ \pi.$  Then $\tilde f_1(\pi(x)-n\alpha)=f(T^{-n} x),$ for all $n$ and $\mu$-almost every $x,$ and $\lim_{j\to\infty} {\tilde g}_{j,1}= \lim_{j\to \infty} \int \tilde f_1\circ R_\alpha^{-n}\, d\eta_j(n).$    By Lemma \ref{weak*limits}, the last limit is the convolution $\tilde f_1*\psi,$ where $\psi:Z\to [0,1]$ is a function satisfying $\int \psi\, dm=d_\nu(A).$   Since $\tilde f_1: Z\to [0,1]$ and $\int \tilde f_1\, dm=\mu(D),$  Lemma \ref{Steinhaus+} now implies that $\tilde g:=\lim_{j\to \infty} \tilde g_{j,1}$ is continuous and the support $U$ of $\tilde g$ has measure at least $\max\{\mu(D),d_\nu(A)\}.$  Then the support of $g:=\lim_{j\to\infty} g_j$ contains $\pi^{-1}(U),$ so we are done. $\square$

\newpage

\section{Examples}\label{examplessection}

Here we find examples of sets $A,B\subset \mathbb Z$ with $d^*(B)>0$ where $A+B$ is not piecewise syndetic.  We will construct these examples from Kronecker systems, via the next lemma and proposition.

\begin{lemma}\label{equivalents}  Let $(Z,\mathscr Z, m, R_\alpha)$ be a Kronecker system, and let $K\subset Z$ be compact.  The following conditions are equivalent.
\begin{enumerate}
\item[(i)] $K$ has nonempty interior.
\item[(ii)] $\{n: n\alpha \in K\}$ is a Bohr set.
\item[(iii)] $\{n: n\alpha \in K\}$ is piecewise syndetic.
\end{enumerate}

\end{lemma}

\textit{Proof.}  (i)$\implies$(ii) follows from the definition of ``Bohr set," and (ii)$\implies$(iii) follows from the fact that Bohr sets are syndetic.

To see that (iii)$\implies$(i), let $K\subset Z$ be compact with $R:=\{n:n\alpha\in K\}$ piecewise syndetic.  Then there exists a finite set $F$ such that $R':=\bigcup_{a\in F} a+R$ is thick.  But $R'$ is the set of return times to a union of translates of $K$: $R'=\{n: n\alpha \in \bigcup_{a\in F} K+a\alpha\}.$   We claim that the thickness of $R'$ implies that $K':=\bigcup_{a\in F} K+a\alpha$ is equal to $Z.$  Since $K'$ is compact, it has open complement.  If $Z\setminus K'$ is nonempty, then $\{n: n\alpha \in Z\setminus K'\}$ is syndetic, which contradicts the fact that $R'$ is thick, so $K'=Z.$  It follows that one of the $K+a\alpha$ has nonempty interior.  Hence $K$ has nonempty interior.  \hfill $\square$

\medskip

\begin{proposition}\label{examples}  Suppose that $A\subset\mathbb Z,$ $\mathbf Z=(Z,\mathscr Z, m,R_\alpha)$ is a Kronecker system, and $\overline{\{n\alpha:n\in A\}}$ has Haar measure $0$ in $Z.$  Then for all $\varepsilon>0,$ there exists $B\subset \mathbb Z$ with $d^*(B)>1-\varepsilon,$ such that $A+B$ is not piecewise syndetic.

Furthermore, if $(\nu_j^{(i)})_{j\in \mathbb N}$ is an equidistributed averaging sequence for each $i\in \mathbb N,$ there exists $B\subset \mathbb Z$ such that $d_{\nu^{(i)}}(B)>1-\varepsilon$ for all $i,$ and $A+B$ is not piecewise syndetic.
\end{proposition}

\textit{Proof.}  Write $E$ for the closure $\overline{\{n\alpha:n\in A\}},$ and by Lemma \ref{Steinhausconverse} let $K\subset Z$ be compact with $m(K)>1-\varepsilon$ such that $E+K$ has empty interior.  By the pointwise ergodic theorem, there exists $z\in Z$ such that
$$
m(K)=\lim_{N\to \infty} \frac{1}{N} \sum_{n=1}^{N} 1_K(z+n\alpha).
$$
Let $B:=\{n: z+n\alpha\in K\},$ so that $d^*(B)\geq m(K)> 1-\varepsilon.$  Then $A+B\subset$ $\{n: z+n\alpha\in E+K\},$ which is not piecewise syndetic, by Lemma \ref{equivalents}.

To prove the second claim, we can, for each $(\nu_j^{(i)})_{j\in \mathbb N},$ pass to a subsequence $(\rho_j^{(i)})_{j\in \mathbb N},$ having the property that the averages $
\lim_{j\to \infty} \int 1_K(z+n\alpha) \, d\rho_j^{(i)}(n)
$
converge to $m(K)$ as $j\to \infty,$ for almost every $z.$  Thus, there is a $z$ that witnesses this convergence for each $i$ simultaneously, and we can proceed as in the previous paragraph, taking $B$ to be $\{n:z+n\alpha\in K\}.$ \hfill $\square$

\medskip

To construct examples via Proposition \ref{examples}, for $\mathbf Z$ we can use the Kronecker system $(b\mathbb Z,b\mathscr Z,m,R_1),$ where $b\mathbb Z$ is the Bohr compactification of $\mathbb Z, b\mathscr Z$ is its Borel $\sigma$-algebra, and $R_1(z)=z+1.$  For $S\subset \mathbb Z,$ let $\tilde S$ be the closure of $S$ in $b\mathbb Z,$  and note that $\tilde S$ is $\overline{\{R_1^n0:n\in S\}}.$

In \cite{DP} the closures $\tilde S_i$ of the following sets $S_i$ in $b\mathbb Z$ are each shown to have Haar measure $0.$

\begin{itemize}
\item  $S_1=$ the set of prime powers (including the primes).
\item  $S_2=$ the set $\{n^2+m^2:n,m\in \mathbb N\}$ of sums of two squares.
\item  $S_3=$ the set of square-full numbers, that is, the set of numbers $n$ so that every exponent in the prime factorization of $n$ is at least two.
\item  $S_4=$ any set of the form $\{\sum_{\varepsilon_i} \varepsilon_i n_i:\varepsilon_i\in \{0,1\}\},$ where $(n_i)_{i\in \mathbb N}$ is a sequence of positive integers satisfying $n_i|n_{i+1}$ for all $i$ and $n_{i+1}/n_i\geq 3$ for all $i.$
\end{itemize}

In \cite{KR}, it was shown that $\tilde S_5$ has Haar measure $0$ whenever $S_5=p(\mathbb Z),$ where $p$ is a polynomial with integer coefficients having degree $2$ or $3.$

By Proposition \ref{examples}, if $A$ is any of the above sets $S_i$ (or the union of finitely many such sets), then there exists a set $B\subset \mathbb Z$ with $d^*(B)>1-\varepsilon$ such that $A+B$ is not piecewise syndetic.  Since the set of square-full numbers includes, for each integer $k\geq 2,$ the set $\mathbb N\text{\textasciicircum} k:=\{n^k: n\in \mathbb N\},$ we may also take $A=\mathbb N\text{\textasciicircum} k$ for $k\geq 2.$

\begin{remark}  In fact, the arguments in \cite{DP} show that closures of the above sets $S_i$, $i\leq 4,$  appropriately embedded in $\prod_{p \text{ prime}} \mathbb Z_p$, have Haar measure $0,$ where $\mathbb Z_p$ is the set of $p$-adic integers, with the usual topology. We could thereby avoid using the Bohr compactification of $\mathbb Z$, which may be desirable given its complexity.
\end{remark}

\section{Questions about sumsets, the Bohr topology, and recurrence}\label{questions}

\subsection*{Sumsets and the Bohr topology}

Recall that $b\mathbb Z$ is the Bohr compactification of $\mathbb Z,$ and $(\nu_j)_{j\in \mathbb N}$ is an equidistributed averaging sequence if and only if $\lim_{j\to \infty} \nu_j=m_{b\mathbb Z}$ in the $\text{weak}^*$ topology of $b\mathbb Z.$

The proof of Theorem \ref{eastheorem} exploited the properties of equidistributed averaging sequences in two different ways.  First, there was Lemma \ref{wm}, which reduced the problem from the setting of a general measure preserving system to the special case of Kronecker systems.  Lemma \ref{wm} can be deduced from the spectral theorem and the following fact,  which is essentially Wiener's lemma:

\medskip

\noindent (W): If $\sigma$ is an atomless probability measure on $\mathbb T$ and $(\nu_j)_{j\in \mathbb N}$ is an equidistributed averaging sequence, then for all $\varepsilon>0, d_\nu\{n: |\int e^{in\theta}\, d\sigma(\theta)|>\varepsilon\}=0.$

\medskip

The second important property of equidistributed averaging sequences is how they project to compact abelian groups.  Lemma \ref{weak*limits} says that when $(\nu_j)_{j\in \mathbb N}$ is an equidistributed averaging sequence and $d_\nu(A)>0,$ the measures $\eta_j:=\nu_j|_A$ are ``large," in the sense that the $\text{weak}^*$-limits of the $\eta_j$ are absolutely continuous with respect to Haar measure on $b\mathbb Z.$  As a consequence we have:

\medskip
\noindent (L): If $(\nu_j)_{j\in \mathbb N}$ is as above, $A\subset \mathbb Z,$ and $d_{\nu}(A)>0,$ then $m_{b\mathbb Z}(\tilde{A})>0,$ where $\tilde{A}$ is the closure of $A$ in $b\mathbb Z.$
 \medskip

In Section \ref{examplessection}, we showed that when the conclusion of (L) fails for $A\subset \mathbb Z,$ there exists $B\subset \mathbb Z$ with $d^*(B)>0$ and $A+B$ is not piecewise syndetic.  We do not understand the situation where the conclusion of  (L) holds but (W) is unavailable.  This situation is not vacuous, for Katznelson (\cite{K}) and Saeki (\cite{Sae}) have produced atomless probability measures $\sigma$ on $\mathbb T$ such that $A_\varepsilon:=\{n:|\int e^{in\theta}\, d\sigma(\theta)|>1-\varepsilon\}$ is dense in $b\mathbb Z$  for all $\varepsilon>0.$  In particular, $d_\nu(A_\varepsilon)=0$ whenever $\nu$ is an equidistributed averaging sequence, while each $\{n\alpha:n\in A_\varepsilon\}$ is dense in $Z$ whenever $\{n\alpha:n\in \mathbb Z\}$ is dense in the compact abelian group $Z.$  We cannot even decide if $A_{1/2}+B$ is piecewise syndetic whenever $d^*(B)>0,$ where the $A_{1/2}$ comes from Katznelson's exmaple $\sigma.$  In general, we ask the following.

\begin{question}\label{bohrquestion}  Let $A\subset \mathbb Z,$ and let $\tilde A$ be the closure of $A$ in $b\mathbb Z.$  Which, if any, of the following implications hold?

\begin{enumerate} \item If $m_{b\mathbb Z}(\tilde A)>0$ and $d^*(B)>0$ then $A+B$ is piecewise syndetic.

\item If $m_{b\mathbb Z}(\tilde A)>0$ and $d^*(B)>0$ then $A+B$ is piecewise Bohr.

    \item If $\tilde A=bZ$ and $d^*(B)>0$ then $A+B$ is thick.
\end{enumerate}
\end{question}

\subsection*{Sets of recurrence}  Call $A\subset \mathbb Z$ a \textit{set of recurrence} if for every measure preserving system $(X,\mathscr X,\mu,T)$ and every $D$ with $\mu(D)>0,$ there exists $n\in A$ with $\mu(D\cap T^{-n}D)>0.$

The following question has been asked in various forms, most recently in Section 9 of \cite{BR}.

\begin{question}\label{denserecurrence} If $A\subset \mathbb Z$ is dense in $b\mathbb Z,$ is $A$ necessarily a set of recurrence?
\end{question}

An affirmative answer to Question \ref{denserecurrence}  would imply an affirmative answer to Part 3 of  Question \ref{bohrquestion}.  The implication is obtained as follows: if every shift of $A$ is a set of recurrence, one can show that whenever $(X,\mathscr X,\mu, T)$ is an ergodic measure preserving system, then $\mu(\bigcup_{a\in A} T^aD)=1$ whenever $\mu(D)>0.$  One can then argue as in the proof of Theorem \ref{eastheorem} to show that $A+B$ is thick whenever $d^*(B)>0.$

\medskip

The next questions might be resolved more easily than Question \ref{bohrquestion}.

\begin{question}\label{thicksums?}  Suppose that $A\subset \mathbb Z$ has the property that $A+B$ is thick whenever $d^*(B)>0$. Must the following be true?

\medskip

$\bullet$ For all ergodic $(X,\mathscr X,\mu,T)$ and all $D\in \mathscr X$ with $\displaystyle \mu(D)>0,$ $\displaystyle \mu\bigl(\bigcup_{a\in A} T^a D\bigr)=1.$
\end{question}

\begin{question}  Suppose that $A\subset \mathbb Z$ has the property that $A+B$ is piecewise syndetic (alternatively, piecewise Bohr) whenever $d^*(B)>0$.  What can be said about $A$?
\end{question}

\subsection*{Two sparse summands.}  Our methods and examples say little about $A~+~B$ when $d^*(A)=d^*(B)=0.$  In particular, let $P$ be the set of primes, and define $d_P(A)=\limsup_{n\to \infty} \frac{|A\cap [1,n]|}{|P\cap [1,n]}$ for $A\subset P.$  We wonder what can be said about $A+B$ when $d_P(A),d_P(B)>0.$  A recent result in \cite{CH} shows that $\bar d(A+A)>0$ whenever $d_P(A)>0.$  Can we conclude that $A+A$ is piecewise syndetic?

\begin{remark}  In \cite{Pav}, R. Pavlov constructs a set $A\subset \mathbb Z$ with $d^*(A)=0$ and the property that $A+B$ is thick whenever $B$ is infinite.  It may be interesting to characterize such $A$ in terms of dynamics.
\end{remark}

\end{document}